\documentclass[bysection, amsart]{standard}

\usepackage{xcolor}
\usepackage{float}
\usepackage{hyperref}
\usepackage{enumitem}

\title{$2$-blocks with defect group of rank $2$}
\author{Pete Gautam}
\address{Department of Mathematics, University of Manchester, Manchester, M13 9PL, United Kingdom}
\email{pratyush.gautam@manchester.ac.uk}


\begin{document}

\begin{abstract}
    We classify all controlled $2$-blocks with defect group of rank $2$ up to Morita equivalence. This provides a proof of Donovan's Conjecture for resistant $2$-generated $2$-groups.
\end{abstract}

\maketitle

\section{Introduction}

Let $p$ be a prime and $(K, \calO, k)$ be a $p$-modular system. Two $k$-algebras $A$ and $B$ are \emph{Morita equivalent} if there is an equivalence between categories of $A$-modules and $B$-modules. A \emph{$p$-block} $B$ of $G$ is an indecomposable factor of $\calO G$. A \emph{defect group} of $B$ is a $p$-subgroup that is maximal among vertices of the modules in $B$, which is unique up to $G$-conjugacy.

A problem in modular representation theory is the classification of $p$-blocks with a fixed defect group up to Morita equivalence. In the case the block has an abelian defect group, there is a classification of $2$-blocks with defect groups that are up to $4$-generated $2$-groups (see \cite{paper:abelian-rank4}). Moreover, there has been some work in cases when the block has a non-abelian defect group, including finite groups of tame representation type (see \cite{paper:quaternion}) and $2$-blocks with a defect group that is a Suzuki $2$-group \cite{paper:eat-suz}. There has also been some work done in the case where the defect group is $2$-generated, including when the defect group is abelian \cite{paper:hom-rank-2} or a minimal non-abelian $2$-group \cite{paper:ekb-block-mna}. In both of these cases, any blocks with such a defect group is controlled. In this paper, we extend the classification by determining all controlled blocks with defect group that is a $2$-generated $2$-group. This work is closely related to Donovan's Conjecture, which was recently claimed to have been proven by OpenAI \cite{OpenAI}.

Let $D$ be a $2$-group of rank $2$. If $\Aut(D)$ is a $2$-group, then any controlled block with defect group $D$ is nilpotent. The structure of nilpotent blocks was determined by Puig in \cite{paper:nilpotent}. The interest is in non-nilpotent, controlled blocks, when $\Aut(D)$ is not a $2$-group. We call such groups \emph{maxperm} $2$-groups of rank $2$. The following is the main result of this paper.
\begin{theorem}
    Let $G$ be a finite group, and let $B$ be a non-nilpotent, controlled block of $G$ with defect group $D$, where $|D| \geq 2^4$ is a maxperm $2$-group of rank $2$. Then $B$ is Morita equivalent to $\calO [D \rtimes C_3]$.
\end{theorem}
\noindent If $D$ is a $2$-generated maxperm $2$-group satisfying $|D| \leq 2^3$, then $D$ is isomorphic to $C_2 \times C_2$ or $Q_8$. The classification of all blocks with such defect groups is already present in the literature (see \cite{paper:li-c2c2} and \cite{paper:ei-q8} respectively).

If $D$ is resistant, then any block with defect group $D$ is controlled. Hence, the theorem has the following corollary.
\begin{corollary}
    Let $G$ be a finite group, and let $B$ be a block of $G$ with defect group $D$, a resistant maxperm $2$-group of rank $2$ with $|D| \geq 2^4$. Then $B$ is Morita equivalent to $\calO D$ or $\calO[D \rtimes C_3]$.
\end{corollary}

% Let $D$ be a $2$-generated $2$-group that is not resistant. Then there may exist a non-controlled block $B$ with defect group $D$. We do not have the classification of blocks in this case.

We briefly explain the structure of the paper. In Section \ref{sec:background}, we recall some background that we use throughout the paper. In Section \ref{sec:maxperm}, we establish elementary properties about maxperm groups, and we study a subfamily of maxperm groups in Section \ref{sec:quasi-ab}. In Section \ref{sec:no-comp}, we deal with the case when the group $G$ has no components, while in Section \ref{sec:comp}, we rule out the case that $G$ has components. The main result is then proven in Section \ref{sec:final}.

\section{Background} \label{sec:background}
\subsection{Group Theory} We make use of the following results throughout the paper without explicit reference.

We recall that $\Phi(G)$ is the intersection of all the maximal subgroups of $G$. If $G$ is a $p$-group, then $G/\Phi(G)$ is elementary abelian.
\begin{lemma}[Coprime Action] \label{lm:coprime}
    Let $G$ and $A$ be finite groups such that $G$ acts coprimely on a $A$, and let $B$ be a $G$-invariant subgroup of $A$. Then the following hold:
    \begin{enumerate}
        \item if $G$ acts trivially on $A/B$ and $B$, then $G$ acts trivially on $A$; 
        \item if $A$ is abelian, then $A = [A, G] \times C_A(G)$; and
        \item if $G$ acts trivially on $A/\Phi(A)$, then $G$ acts trivially on $A$.
    \end{enumerate}
\end{lemma}
\begin{proof}
    See \cite[Section 8]{ks06}.
\end{proof}

We can strengthen (3) in Lemma \ref{lm:coprime} by the following result of Burnside.

\begin{lemma}[Burnside] \label{lm:burnside}
    Let $D$ be a finite $p$-group. Fix $A \leq \Aut(D)$ to be the subgroup consisting of all $\theta \in \Aut(D)$ that act trivially on $D/\Phi(D)$. Then $A$ is a normal $p$-subgroup of $\Aut(D)$.
\end{lemma}
\begin{proof}
    See \cite[Theorem I.5.1.4]{gor07}.
\end{proof}
\noindent We will denote the subgroup $A$ described in Lemma \ref{lm:burnside} as $C_{\Aut(D)}(D/\Phi(D))$.

Let $G$ be a finite group. We say that $A \leq G$ is a \emph{component} of $G$ if $A$ is a subnormal quasisimple subgroup of $G$. We define the \emph{layer} $E(G)$ to be the subgroup generated by the components of $G$. The \emph{Fitting subgroup} of $G$ is the largest nilpotent normal subgroup of $G$. It is the direct product of $O_p(G)$ for every prime $p$ dividing $|G|$. The \emph{generalized Fitting subgroup} of $G$ is defined by $F^*(G) = F(G) E(G)$. It is well-known that $C_G(F^*(G)) \leq F^*(G)$ (see, for example, \cite[31.13]{book:aschbacher}).

We recall that the derived subgroup $D'$ is the smallest subgroup $A \leq D$ such that $D/A$ is abelian. Moreover, we say that a finite group $G$ is \emph{perfect} if $G = G'$. Every quotient of a perfect group is perfect.

\begin{lemma} \label{lm:perfect-cen}
    Let $Z \leq Z(G)$ be such that $G/Z$ is a perfect group, and set $G_0 := G'$. Then $G_0$ is a perfect group with $G = G_0Z$.
\end{lemma}
\begin{proof}
    See \cite[1.5.3]{ks06}.
\end{proof}

\begin{lemma}[Gr{\"u}n's Lemma] \label{lm:grun}
    Let $G$ be a perfect group. Then $Z(G/Z(G)) = 1$.
\end{lemma}
\begin{proof}
    See \cite[Satz 4]{paper:grun}.
\end{proof}

\begin{lemma} \label{lm:rank2-max-rank3}
    Let $D$ be a $2$-group of rank $2$. If $D$ has an elementary abelian subgroup $A$ of index $2$, then $|D| \leq 2^4$.
\end{lemma}
\begin{proof}
    Let $x \in D \setminus A$ and $y \in A \setminus \Phi(D)$. Since $D$ is $2$-generated, it follows that $D = \langle x, y\rangle$. Set $z := [x, y]$. Then $z$ lies inside $\Phi(D) \leq A$, meaning that $z^2 = 1$. Moreover, since $A \leq C_D(z)$, we find that $z$ has at most $2$ conjugates in $D$, with $z^D \subseteq A$. Thus,     \[|\langle z^D \rangle| \leq 2^2.\]
    As $D$ is generated by $x$ and $y$, we deduce that $D/\langle z^D \rangle$ is abelian. Since $D'$ is a normal subgroup of $D$ containing $z$, we infer that $D' = \langle z^D \rangle$. 
    
    Now $D/D'$ is an abelian group of rank $2$ with $A/D'$ elementary abelian of index $2$. This forces 
    \[|D : D'| \leq 2^3.\]
    Combining the two equations, we deduce that $|D| \leq 2^5$. We make use of GAP \cite{GAP4} to infer that there is no group of order $2^5$ satisfying this property and conclude that $|D| \leq 2^4$.
\end{proof}

\subsection{Fusion Systems}
See \cite{book:cra} or \cite{ako} for background on fusion systems. When working with a fusion system $\calF$, we will try to understand its structure by understanding the subgroups that are \emph{essential} in $\calF$, as a consequence of the Alperin-Goldschmidt Theorem.
\begin{theorem}[Alperin-Goldschmidt]
    Let $\calF$ be a saturated fusion system on a finite $p$-group $D$. Then there exist a set of \emph{essential subgroups}, denoted by $\calE(\calF)$, such that
    \[\calF = \langle \Aut_\calF(D), \Aut_\calF(E) \mid E \in \calE(\calF) \rangle.\]
\end{theorem}
\begin{proof}
    See \cite[Theorem I.3.5]{ako}.
\end{proof}
\noindent We say that a fusion system $\calF$ is \emph{controlled} if $\calE(\calF) = \varnothing$, and that $\calF$ is \emph{nilpotent} if $\calF$ is controlled and $\Aut_\calF(D) = \Inn(D)$. A nilpotent fusion system is equal to $\calF_D(D)$. A finite $p$-group $D$ is called \emph{resistant} if for every saturated fusion system $\calF$ on $D$, $\calF$ is controlled.

Let $\calF$ be a fusion system on $D$, and $A \leq D$. We define the \emph{centralizer fusion system} $C_\calF(A)$ to be a fusion system on $C_D(A)$, where a map $\phi \colon X \to Y$ in $\calF$ lies in $C_\calF(A)$ if there exists an extension $\hat{\phi} \colon AX \to AY$ of $\phi$ in $\calF$ such that $\hat{\phi}|_A = 1_A$. If $A$ is fully $\calF$-centralized, then $C_\calF(A)$ is saturated by \cite[Theorem 4.28]{book:cra}.

\subsection{Block Theory} See \cite{book:alperin-local-rep-theory} or \cite{book:navarro} for background on block theory. A block $B$ with defect group $D$ has a saturated fusion system $\calF$ (see, for example, \cite[Part IV]{ako}). If $D$ is a defect group of a block $B$ of some group $G$, then there is a corresponding (saturated) fusion system $\calF := \calF_{(D, b)}(G, B)$ on $D$, where $b$ is a block of $DC_G(D)$ with $b^G = B$. We say that $B$ is \emph{controlled} if $\calF$ is controlled. Similarly, we say that $B$ is \emph{nilpotent} if $\calF$ is nilpotent.

\begin{lemma} \label{lm:idxp-nilpotent}
    Let $G$ be a group, $B$ be a $p$-block of $G$ and let $N \normalIn G$ be such that $G/N$ is a $p$-group. Let $b$ be a block of $N$ covered by $B$. Then $B$ is nilpotent if and only if $b$ is nilpotent.
\end{lemma}
\begin{proof}
    See \cite[Proposition 6.5]{paper:kulshammer-puig}.
\end{proof}

\begin{lemma} \label{lm:cen-prod-nilpotent}
    Let $G$ be a finite group that is a central product of groups $G_i \normalIn G$ for $1 \leq i \leq n$. Let $B$ be a block of $G$, and let $b_i$ be the (unique) block of $G_i$ covered by $B$. Then $B$ is nilpotent if and only if each $b_i$ is nilpotent.
\end{lemma}
\begin{proof}
    See \cite[Lemma 7]{paper:ekb-block-mna}.
\end{proof}

\begin{lemma}
    Let $G$ be a finite group, and let $H \normalIn G$. Let $B$ be a block of $G$, and let $b$ be a block of $H$ covered by $B$. If $D$ is a defect group of $B$, then $D \cap H$ is $G$-conjugate to a defect group of $b$.
\end{lemma}
\begin{proof}
    See \cite[Theorem 15.1]{book:alperin-local-rep-theory}.
\end{proof}

\begin{proposition} \label{prp:block-domination}
    Let $G$ be a finite group, and let $N \normalIn G$ be a $p$-group. Write $\overline{G} := G/N$. If $B$ is a $p$-block of $G$ with defect group $D$, then there exists a block $\overline{B}$ of $\overline{G}$ with defect group $\overline{D}$. Moreover if $N \leq Z(G)$, then $B$ is nilpotent if and only if $\overline{B}$ is nilpotent.
\end{proposition}
\begin{proof}
    See \cite[Theorem 9.9 (b)]{book:navarro} and \cite[Lemma 2]{paper:nilpotent-dom}.
\end{proof}

% For a group $G$, the principal block of $kG$ is denoted by $B_0(G)$.
% \begin{theorem}[Brauer's Third Main Theorem]
%     Let $G$ be a finite group, and let $b$ be a block of $kH$, where $H \leq G$ such that $b^G$ is defined. Then $b^G$ is the principal block of $kG$ if and only if $b$ is the principal block of $kH$.
% \end{theorem}
% \begin{proof}
%     See \cite[Theorem 6.7]{book:navarro} and \cite[Theorem IV.5.9]{ako}.
% \end{proof}
% \noindent As a corollary of this theorem, we deduce that if $B_0$ is the principal block of $G$, and $\calF$ the block fusion system of $B_0$, then $\calF = \calF_D(G)$, where $D$ is a Sylow $p$-subgroup of $G$. This allows us to conclude, for example, when principal blocks are not nilpotent.

% \begin{corollary}
%     The principal block of $\PSL_2(3)$ is not nilpotent.
% \end{corollary}

We will mainly be working with blocks that obey the following properties.
\begin{proposition} \label{prp:reduced}
    Let $G$ be a finite group and let $B$ a block of $\calO G$ with defect group $D$. Then there exists a finite group $H$ and a block $C$ of $\calO H$ such that $B$ is Morita equivalent to $C$, a defect group $P_H$ of $C$ is isomorphic to $D$, and for any $N \normalIn G$,
    \begin{enumerate}[label=(R\arabic*)]
        \item $C$ covers a unique block of $N$; and
        \item if $C$ covers a nilpotent block of $\calO N$, then $N \leq O_p(H) Z(H)$ and $O_{p'}(H) \leq H'$.
    \end{enumerate}
\end{proposition}
\begin{proof}
    See \cite[Proposition 6.1]{paper:reduced-blocks}.
\end{proof}
\noindent The pair $(H, C)$, or simply $C$, satisfying the conditions (R1) and (R2) above is called \emph{reduced}.

\section{Maxperm Groups} \label{sec:maxperm}
Let $D$ be a $p$-group. We define the \emph{rank} $r \geq 1$ of $D$ to be the value satisfying $|D : \Phi(D)| = p^r$. If $D$ has rank $r$, then any minimal generating set for $D$ contains $r$ elements.

Since $\Phi(D)$ is the intersection of all maximal subgroups of $D$, there is a one-to-one correspondence between the maximal subgroups of $D$ and the maximal subgroups of $D/\Phi(D)$. Thus, $D$ is a $p$-group of rank $r$ if and only if $D$ has $\frac{p^r-1}{p-1}$ maximal subgroups. In particular, $D$ is a $2$-group of rank $2$ if and only if $D$ has $3$ maximal subgroups.

\begin{definition}
    Let $D$ be a $p$-group. We say that $D$ is a \emph{maxperm group} if there exists a $\theta \in \Aut(D)$ that permutes the maximal subgroups of $D$ transitively. If $D$ is a maxperm group, then we also say that $D$ is a \emph{maxperm group with respect to $\theta$}.
\end{definition}
Since the maximal subgroups of $D$ are $\Aut(D)$-conjugate, the following result follows immediately.
\begin{lemma} \label{lm:maxperm-max-iso}
    Let $D$ be a maxperm $p$-group. Then the maximal subgroups of $D$ are all isomorphic. Moreover, any maximal subgroup of $D$ is not characteristic in $D$.
\end{lemma}

Let $D$ be a maxperm $p$-group of rank $r$. The quotient
\[O := \Aut(D)/C_{\Aut(D)}(D/\Phi(D))\]
acts faithfully on $D/\Phi(D)$. If $D$ has rank $r$, then $D/\Phi(D)$ is elementary abelian of order $p^r$. This implies that $O$ embeds into $\GL_r(p)$. Since $\theta$ permutes the maximal subgroups of $D$ transitively, the order of $\theta$ is a multiple of $\frac{p^r-1}{p-1}$. In particular, if $p=2$, $\theta$ must be a Singer cycle as an element of $\GL_r(p)$. 

We now establish some elementary properties about these groups that we require.
\begin{lemma} \label{lm:maxperm-char-in-frat}
    Let $D$ be a maxperm $p$-group of rank $r$ with respect to $\theta \in \Aut(D)$. If $N$ is a proper subgroup of $D$ that is invariant under $\langle \theta \rangle$, then $N \leq \Phi(D)$.
\end{lemma}
\begin{proof}
    Assume that $N < D$ is invariant under $\langle \theta \rangle$ of $D$. Since $N \neq D$, there exists a maximal subgroup $M$ such that $N \leq M$. Then we find that 
    \[N = \theta^i(N) \leq \theta^i(M)\]
    for $1 \leq i \leq p^r - 1$. Hence, $N$ is contained in every maximal subgroup of $D$. This shows that $N \leq \Phi(D)$.
\end{proof}

\begin{lemma} \label{lm:maxperm-quo}
    Let $D$ be a maxperm $p$-group of rank $r$ with respect to $\theta \in \Aut(D)$, and let $N \normalIn D$ be a proper subgroup that is invariant under $\langle \theta \rangle$. Then $D/N$ is a maxperm $p$-group of rank $r$ with respect to $\theta$. In particular, if $D/N$ is abelian, then it is a homocyclic $p$-group of rank $r$.
\end{lemma}
\begin{proof}
    We know by Lemma \ref{lm:maxperm-char-in-frat} that $N \leq \Phi(D)$. As such, there is a one-to-one correspondence between the maximal subgroups of $D$ and the maximal subgroups of $D/N$. In particular, $D/N$ has $p^r-1$ maximal subgroups. We deduce that $D/N$ has rank $r$. By the Third Isomorphism Theorem, we infer that $\Phi(D/N) = \Phi(D)/N$.
    
    Let $M$ be a maximal subgroup of $D$. Since $N \leq \Phi(D) \leq M$, we infer that
    \[\theta^i(M/N) = \theta^i(M)/N\]
    for $1 \leq i \leq p^r - 1$. Thus, $\theta$ must permute the maximal subgroups of $D/N$ transitively. We conclude that $D/N$ is a maxperm $p$-group of rank $r$ with respect to $\theta$.

    Finally, assume that $D/N$ is abelian. We have shown that $D/N$ is a maxperm group of rank $r$. In particular, the maximal subgroups of $D/N$ are all isomorphic. This forces $D/N$ is a homocyclic $p$-group of rank $r$.
\end{proof}

In this paper, we be considering the maxperm groups of rank $2$. The following result motivates our study.
\begin{lemma} \label{lm:maxperm-rank2}
    Let $D$ be a $2$-group of rank $2$. If $\Aut(D)$ is not a $2$-group, then $D$ is maxperm.
\end{lemma}
\begin{proof}
    Since $D$ is a $2$-group of rank $2$, $D/\Phi(D)$ is an elementary abelian group of order $4$. The subgroup
    \[\Aut(D)/C_{\Aut(D)}(D/\Phi(D))\]
    acts faithfully on $D/\Phi(D)$. This implies that the subgroup embeds into $\GL_2(2) \cong S_3$. By Lemma \ref{lm:burnside}, we further know that $C_{\Aut(D)}(D/\Phi(D))$ is a $2$-group. As such, if $\Aut(D)$ is not a $2$-group, there exists a $3$-element $\theta \in \Aut(D)$. By construction, $\theta$ acts non-trivially on $D/\Phi(D)$. Hence, $\theta$ permutes the maximal subgroups of $D/\Phi(D)$. Taking preimage in $D$, we deduce that $\theta$ permutes the three maximal subgroups of $D$. Thus, $D$ is a maxperm group with respect to $\theta$.
\end{proof}

A $p$-group $D$ is \emph{minimal non-abelian} if every maximal subgroup of $D$ is abelian, but $D$ is non-abelian. The following lemma characterizes when a minimal non-abelian group is maxperm. We recall that if $D$ is a group with a central subgroup $Z$ such that $D/Z$ is cyclic, then $D$ is abelian.
\begin{lemma}\label{lm:maxperm-small}
    Let $D$ be a non-abelian maxperm $2$-group of rank $2$ with respect to $\theta$. Then $D$ is minimal non-abelian if and only if there exists a subgroup $Z \leq Z(D)$ invariant under $\langle \theta \rangle$ such that $|D : Z| \leq 2^3$.
\end{lemma}
\begin{proof}
    Assume first that $|D : Z| \leq 2^3$. By Lemma \ref{lm:maxperm-quo}, we know that $D/Z$ is a maxperm $2$-group of rank $2$ of order $\leq 2^3$. This implies that $D/Z$ is isomorphic to either $C_2^2$ or $Q_8$. In particular, every maximal subgroup $M$ of $D$ is such that $M/Z$ is cyclic. This means that $M$ is abelian, and we conclude that $D$ is minimal non-abelian.

    Assume now that $D$ is minimal non-abelian. Consider the subgroup $C := C_D(\Phi(D))$. The subgroup $\Phi(D)$ is contained in each maximal subgroup of $D$, which are all abelian. We deduce that $C$ contains every maximal subgroup of $D$, so that $C = D$. This implies that $\Phi(D) \leq Z(D)$. Hence, we take $Z := \Phi(D)$, which satisfies $|D : Z| = 2^2$.
\end{proof}

$2$-blocks with defect groups that are both maxperm and minimal non-abelian were considered by Eaton, K{\"u}lshammer and Sambale in \cite{paper:ekb-block-mna}. A pivotal lemma they use to prove this result is \cite[Lemma 5]{paper:ekb-block-mna}. The following extends their result for maxperm $2$-groups.

\begin{lemma} \label{lm:smaller-not-nilpotent}
    Let $G$ be a finite group, $B$ be a non-nilpotent, controlled block of $G$ with defect group $D$ a maxperm $2$-group of rank $2$. Let $H \normalIn G$ and $b$ be a block of $H$ covered by $B$ with defect group $D_0$. If $|D_0| < |D|$, then $b$ must be nilpotent.
\end{lemma}
\begin{proof}
    Without loss of generality, we assume that $D_0 = D \cap H$. Since $H$ is normal in $G$, it follows that $D_0$ is normal in $N_G(D)$. We see that $B$ is a non-nilpotent, controlled block with inertial quotient $\langle \theta \rangle$, where $\theta \in \Aut_G(D)$ is a $3$-element. By assumption, $D_0$ is a proper subgroup of $D$ that is invariant under $\langle \theta \rangle$. Applying Lemma \ref{lm:maxperm-char-in-frat}, we deduce that $D_0 \leq \Phi(D)$. In that case, we have
    \[[D, N_H(D)] \leq D \cap H = D_0 \leq \Phi(D).\]
    As such, $N_H(D)$ acts trivially on $D/\Phi(D)$. Using Lemma \ref{lm:burnside}, we infer that $N_H(D)/C_H(D)$ is a $2$-group. 
    
    Now let $\beta$ be the unique block of $DH$ that covers $B$. We may assume that $D$ is a defect group of $\beta$ by \cite[Theorem E]{paper:p-solvable}. We next let $\beta_D$ be a block of $DC_{DH}(D)$ such that $(\beta_D)^{DH} = \beta$. By the previous paragraph, it follows that $N_H(D, \beta_D)/C_H(D)$ and $N_{DH}(D, \beta_D)/C_{DH}(D)$ are  $2$-groups. As such, $\beta$ has trivial inertial quotient. Thus, the block $\beta$ must be nilpotent. Since $DH/H$ is a $2$-group, we conclude by Lemma \ref{lm:idxp-nilpotent} that $b$ is indeed nilpotent.
\end{proof}

When proving the main result, we will need to identify all maxperm $2$-groups of rank $2$ that can be a Sylow $2$-subgroup of a simple group. To do so, we make use of the following analysis of quasisimple groups that have a Sylow $2$-subgroup that is $2$-generated. This is based on the work by Navarro et al in \cite{paper:rank2-chars}.
\begin{table}[H] 
    \centering
    \begin{tabular}{c|c|c|c}
        $G$ & $q$ & Structure of $D$ & Result \\
        \hline
        $A_7$ & $-$ & $D_8$ & Proposition 4.2 \\
        ${^2 \mathrm{F}_4(2)'}$ & $-$ & $[2^{11}]$ & Lemma 4.3  \\
        $\PSL_2(q)$ & $q \equiv \pm 1 \bmod{8}$ & $D_8$ & Proposition 4.12 \\
        $\PSL_2(q)$ & $q \equiv \pm 3 \bmod{8}$ & $C_2^2$ & Proposition 4.12 \\
        $\PSL_3(q)$ & $q \equiv 1 \bmod{4}$ & $C_{2^a} \wr C_2$ & Proposition 4.12 \\
        $\PSL_3(q)$ & $q \equiv 3 \bmod{4}$ & Semi-dihedral & Proposition 4.12 \\
        $\PSU_3(q)$ & $q \equiv 1 \bmod{4}$ & Semi-dihedral & Proposition 4.12 \\
        $\PSU_3(q)$ & $q \equiv 3 \bmod{4}$ & $C_{2^a} \wr C_2$ & Proposition 4.12 
    \end{tabular}
    \caption{All simple groups $G$ with $D\in \Syl_2(G)$ of rank $2$ (up to isomorphism). A reference to the specific result identifying $G$ in \cite{paper:rank2-chars} is also provided, along with a description of the structure of $D$. The value $q$ is a prime power, and $a \geq 2$ depends on $q$. The structure of $D$ for the classical groups of dimension $\leq 3$ is given in \cite[Chapter 6.5]{book:gls3}.}
    \label{tbl:rank2}
\end{table}

\begin{lemma} \label{lm:sylow-choices}
    Let $D$ be a maxperm $2$-group of rank $2$, and let $G$ be a simple group with $D \in \Syl_2(G)$. Then $G \cong \PSL_2(q)$ with $q \equiv \pm 3 \bmod{8}$.
\end{lemma}
\begin{proof}
    The $2$-groups of rank $2$ that are isomorphic to a Sylow $2$-subgroup of a simple group are given in Table \ref{tbl:rank2}. In most cases, the following arguments contradict that $D$ is a maxperm $2$-group:
    \begin{itemize}
        \item if $D \cong D_8$, then $\Aut(D) \cong D_8$ is a $2$-group;
        \item if $D$ is a Sylow $2$-subgroup of ${^2 \mathrm{F}_4(2)'}$, then we make use of MAGMA \cite{magma} to compute that $\Aut(D)$ is a $2$-group;
        \item if $D$ is semi-dihedral, then we know, for example, by \cite[Exercise 2.1 (9)]{book:primepower} that the maximal subgroups of $D$ are not isomorphic, contradicting Lemma \ref{lm:maxperm-max-iso};
        \item if $D \cong C_{2^a} \wr C_2$ for $a \geq 2$, then $D$ has a maximal subgroup that is characteristic by \cite[Theorem 9.12]{paper:wreath-products}, contradicting Lemma \ref{lm:maxperm-max-iso}.
    \end{itemize}
    We are left with only one choice -- $G \cong \PSL_2(q)$ for $q \equiv \pm 3 \bmod{8}$, as desired.
\end{proof}

% Let $D$ be a minimal non-abelian group. We know by \cite[Theorem 2.4]{paper:sambale} that any such $D$ is resistant. The argument in \cite{paper:ekb-block-mna} makes use of \cite[Lemma 5]{paper:ekb-block-mna} to show that $D$ must be a defect group of a block of a quasisimple group $G$. They then make use of the classification of finite simple groups to deduce that no such $(G,B)$ exists.

% Now let $D$ be a maxperm $2$-group of rank $2$, and assume that $D$ is resistant. In this case, we can use a similar argument to show that $D$ must be a defect group of a block $B$ a quasisimple group $G$ (see Lemma \ref{lm:resistant-no-comp}). Then the block $B$ is controlled, so we can make use of \cite{paper:j-an-controlled-2block} to conclude that no such $(G,B)$ exists. 

Let $D$ be a maxperm $2$-group of rank $2$. If $D$ is abelian, then Lemma \ref{lm:maxperm-quo} reveals that $D$ is homocylic group of rank $2$. $2$-blocks with defect group homocyclic group of rank $2$ have been classified in \cite{paper:hom-rank-2}. As such, we shall only be considering blocks with non-abelian defect groups in this paper.

\section{Quasi-abelian groups} \label{sec:quasi-ab}
In this section, we take a close look at a subfamily of maxperm $2$-groups of rank $2$. We shall make use of this in Section \ref{sec:comp}.

\begin{definition}
    Let $D$ be a $p$-group. We say that $D$ is \emph{quasi-abelian} if $D$ has rank $2$ and nilpotency class $2$.
\end{definition}

\begin{lemma}
    Let $D$ be a quasi-abelian group. Then $D'$ is cyclic.
\end{lemma}
\begin{proof}
    Let $x, y \in D$ generate $D$. Consider the subgroup $Z := \langle [x, y] \rangle$. Since $D$ has nilpotency class $2$, $Z \leq Z(D)$ is a normal subgroup of $D$. Moreover, $D/Z$ is abelian by construction. We deduce that $D'$ is cyclic.
\end{proof}

\begin{lemma} \label{lm:minab-qab}
    Let $D$ be a minimal non-abelian group. Then $D$ is quasi-abelian.
\end{lemma}
\begin{proof}
    In \cite[Lemma 12.1]{book:sambale}, Sambale shows that a minimal non-abelian group has rank $2$ and $D'$ has order $p$. Since $D'$ is a normal subgroup, we further know that $D' \leq Z(D)$. As such, $D$ has nilpotency class $2$, meaning that $D$ is quasi-abelian.
\end{proof}
\noindent We think of quasi-abelian groups as a generalisation of minimal non-abelian groups, where instead of requiring that $|D'| = p$, we only need that $D$ has nilpotency class $2$.

\begin{lemma} \label{lm:qab-small}
    Let $D$ be a maxperm $2$-group of rank $2$. If $|D| \leq 2^6$, then $D$ is quasi-abelian.
\end{lemma}
\begin{proof}
    We verify this result using GAP \cite{GAP4}.
\end{proof}

\begin{lemma}
    Let $D$ be a quasi-abelian group, and let $A \leq D$. Then $A$ has rank $\leq 3$.
\end{lemma}
\begin{proof}
    Let $A \leq D$. Since $D/D'$ is an abelian group of rank $2$, it follows that $AD'/D'$ has rank $\leq 2$. As such, there exist $x, y \in A$ such that 
    \[A / (A \cap D') \cong AD' / D' = \langle xD', yD' \rangle.\]
    Next, take a generator $z \in A \cap D'$. Then we find that $A = \langle x, y, z \rangle$, as desired.
\end{proof}

\begin{proposition} \label{prp:metabelian-res}
    Let $D$ be a maxperm $2$-group of rank $2$. If $D$ is quasi-abelian, then $D$ is resistant.
\end{proposition}
\begin{proof}
    Let $\calF$ be a saturated fusion system on $D$. We show that $\calE(\calF) = \varnothing$. Assume that there exists an $E \in \calE(\calF)$. Since $D$ has nilpotency class $2$, it follows that $E \normalIn D$. But $E$ also has rank $\leq 3$, so we are forced to have that $E$ is a maximal subgroup of $D$ by \cite[Lemma 3.3]{paper:oliver-srank4}.
    
    Set $A := E'$. We consider $\calF_0 := C_\calF(A)$, which is a saturated subsystem of $\calF$ since $A \leq Z(D)$. Let $\theta \in \Aut(D)$ be a $3$-element. Since $D'$ is a cyclic $2$-group, we find that $\theta$ centralizes $D'$. As such, $\theta$ also centralizes $A$. If $\psi \in \Aut_\calF(D)$ is a $3$-element, this argument shows that $\psi$ lies inside $\calF_0$. We deduce that $\Aut_\calF(D) = \Aut_{\calF_0}(D)$.
    
    Let $M$ be a maximal subgroup of $D$ such that $M \in \calE(\calF)$. Then $M$ and $E$ are conjugate in $\langle \theta \rangle$. Since $\theta$ centralizes $A = E'$, we infer that $A = M'$ as well. By \cite[Lemma 3.3]{paper:oliver-srank4}, we know that $\Out_\calF(M) \cong S_3$. Take $r \in \Aut_\calF(M)$ a $3$-element. Since $r$ normalizes $A$, which is a cyclic $2$-group, we deduce again that $r$ lies inside $\calF_0$. This implies that $\Aut_\calF(M) = \Aut_{\calF_0}(M)$. Since every essential subgroup in $\calF$ is a maximal subgroup of $D$, we deduce that
    \[\calF = \langle \Aut_\calF(D), \Aut_\calF(M) \mid M \in \calE(\calF) \rangle = \calF_0 = C_\calF(A)\]
    by Alperin-Goldschmidt. 

    Let $\overline{\calF} := \calF/A$ be the quotient fusion system on $\overline{D} := D/A$. Then \cite[Theorem 5.6]{book:cra} reveals that $\overline{E} \in \calE(\overline{\calF})$. This implies that $\overline{D}$ cannot be abelian. By Lemma \ref{lm:maxperm-quo}, $\overline{D}$ is a maxperm group. Moreover, every maximal subgroup of $\overline{D}$ is abelian by construction. Thus, $\overline{D}$ is a minimal non-abelian group. We know by Lemma \ref{lm:maxperm-quo} that $\overline{D}$ is also maxperm. But then \cite[Theorem 2.4]{paper:sambale} tells us that $\overline{D}$ is resistant, contradicting that $\overline{E} \in \calE(\overline{\calF})$.
\end{proof}

We finish this section with a different characterization of quasi-abelian groups.
\begin{lemma} \label{lm:quasiab-equiv}
    Let $D$ be a non-abelian maxperm $2$-group of rank $2$. Then $D$ is quasi-abelian if and only if there exist $A, T \leq D$ such that:
    \begin{enumerate}
        \item $A \normalIn D$;
        \item $D = AT$; 
        \item $A$ and $TZ(D)/Z(D)$ are abelian.
    \end{enumerate}
\end{lemma}
\begin{proof}
    Assume first that $D$ is quasi-abelian. Since $D$ has nilpotency class $2$, it follows that $D' \leq Z(D)$. Then Lemma \ref{lm:maxperm-quo} reveals that $D/D'$ is a homocyclic group of rank $2$. Let $A$ and $T$ be subgroups of $D$ containing $D'$ such that $D/D'$ is the direct product of the cyclic groups $A/D'$ and $T/D'$. Since $D' \leq Z(D)$ and $A/D'$ is cyclic, we deduce that $A$ is abelian. Similarly, $T$, and so $TZ(D)/Z(D)$, is abelian. We further have that $A \normalIn D$, meaning that $D$ satisfies conditions (1) -- (3).
    
    We now prove the converse. Since $D$ is known to have rank $2$, it suffices to show that $D$ has nilpotency class $2$. We consider the quotient $D/A$, which is either cyclic or has rank $2$.

    Assume first that $D/A$ has rank $2$. Then the three maximal subgroups of $D$ all contain $A$. This implies that $A \leq \Phi(D)$. Since $D = AT$, we deduce that $D = T$. By assumption, $TZ(D)/Z(D)$ is abelian, meaning that $D = T$ has nilpotency class $2$.
    
    Otherwise, $D/A$ is cyclic. In that case, $A$ is contained in precisely one maximal subgroup $M$ of $D$. This implies that $A \nleq \Phi(D)$. If $\theta \in \Aut(S)$ is a $3$-element, let $A_0$ be the intersection of the three conjugates of $AZ(D)$ under $\langle \theta \rangle$. We consider $C := C_D(A_0)$. As $AZ(D)$ is abelian, we know that $A \leq AZ(D) \leq C$. In particular, $C \nleq \Phi(D)$. By construction, $A_0$ is invariant under $\langle \theta \rangle$. We deduce by Lemma \ref{lm:maxperm-char-in-frat} that $C = D$. Thus, $A_0 \leq Z(D)$. Moreover, $A_0$ is the intersection of the three conjugates of $AZ(D)$ by $\langle \theta \rangle$, meaning that $A_0 = Z(D)$. Since $D/A$ is cyclic, we have that $D' \leq A \leq AZ(D)$. As $D'$ is invariant under $\langle \theta \rangle$, we conclude that $D' \leq A_0 = Z(D)$. Hence, $D$ has nilpotency class $2$.
\end{proof}

\section{No components} \label{sec:no-comp}
In this section, we prove the following result.

\begin{proposition} \label{prp:no-comp}
    Let $G$ be a finite group and $(G, B)$ be a reduced controlled block. Assume that $B$ has defect group $D$ that is a non-abelian maxperm $2$-group of rank $2$, with $O_2(G) < D$. Then $E(G) \neq 1$.
\end{proposition}
By Proposition \ref{prp:reduced}, we know that every block $B$ is Morita equivalent to a reduced block. If $O_2(G) = D$, then $G$ has normal defect. Since $N_G(D)/DC_G(D) \cong C_3$, we deduce by \cite{paper:kul-normals} that $B$ is Morita equivalent to $\calO [D \rtimes C_3]$. As such, it suffices to consider the case when $O_2(G) < D$. Since $(G, B)$ is reduced, $B$ cannot be nilpotent. Since $D$ is a maxperm group, it follows that $N_G(D)/C_G(D)$ has a $3$-element.

We prove Proposition \ref{prp:no-comp} in several steps. We assume the setup in Proposition \ref{prp:no-comp} throughout the section. We further set $\overline{G} := G/F^*(G)$.

\begin{lemma}
    $O_2(\overline{G}) = 1$.
\end{lemma}
\begin{proof}
    Let $A$ be the full preimage of $O_2(\overline{G})$ in $G$, and fix $P \in \Syl_2(A)$. Since $O_{2'}(A) = O_{2'}(G) \leq Z(G)$, it follows that $A = P O_{2'}(A)$ is a normal nilpotent subgroup of $G$. As $B$ is reduced, it follows that $A = F(G) = F^*(G)$, and so $O_2(\overline{G}) = 1$.
\end{proof}

\begin{lemma}
    $O_{2'}(\overline{G}) = 1$.
\end{lemma}
\begin{proof}
    Let $A$ be the full preimage of $O_{2'}(\overline{G})$. Let $b$ be the unique block of $A$ covered by $B$. Then $b$ has defect group $O_2(G) < D$. Applying Lemma \ref{lm:smaller-not-nilpotent}, we find that $b$ is nilpotent. Since $B$ is reduced, this forces $A \leq F^*(G)$. In particular, we have that $O_{2'}(\overline{G}) = 1$.
\end{proof}

\begin{lemma}
    $E(\overline{G}) \neq 1$.
\end{lemma}
\begin{proof}
    We have found that
    \[F(\overline{G}) = O_2(\overline{G}) O_{2'}(\overline{G}) = 1.\]
    If $E(\overline{G}) = 1$ as well, then $F^*(\overline{G}) = F(\overline{G}) E(\overline{G}) = 1$, and so
    \[\overline{G} = C_G(F^*(\overline{G})) \leq F^*(\overline{G}) = 1.\]
    Thus, $G = F^*(G) = O_2(G) Z(G)$. But then $D = O_2(G)$, a contradiction.
\end{proof}

\begin{lemma}
    $\overline{D} \leq E(\overline{G})$.
\end{lemma}
\begin{proof}
    Let $E$ be the full preimage of $E(\overline{G})$ in $G$. If $D$ is not contained inside $E$, then by Lemma \ref{lm:smaller-not-nilpotent}, it follows that the block $B_E$ of $E$ covered by $B$ is nilpotent. But then $E = F^*(G)$, contradicting that $E(\overline{G}) \neq 1$.
\end{proof}

\begin{lemma}
    $B$ is of maximal defect.
\end{lemma}
\begin{proof}
    Let $b$ be the unique block of $F^*(G)$ covered by $B$. If $P := O_2(G)$, then $H := PC_G(P) = F^*(G)$. By \cite[Corollary 9.21]{book:navarro}, we find that $B = b^G$ is the unique block of $G$ covering $b$. Since a defect group of $b$ is $O_2(G)$, which is a Sylow $2$-subgroup of $F^*(G)$, we conclude by \cite[Theorem 9.17]{book:navarro} that $B$ also has maximal defect.
\end{proof}

\begin{lemma}
    $E(\overline{G}) \cong \PSL_2(q)$, with $q \equiv \pm 3 \bmod{8}$.
\end{lemma}
\begin{proof}
    Since $O_2(G)$ is invariant under $N_G(D)$, Lemma \ref{lm:maxperm-char-in-frat} reveals that $\overline{D}$ is a maxperm group of rank $2$. As $F(\overline{G}) = 1$, $E(\overline{G})$ is a direct product of simple groups. A Sylow $2$-subgroup of a simple group cannot be cyclic. This means that if $\overline{G}$ has multiple components, then a Sylow $2$-subgroup of $E(\overline{G})$ has rank $\geq 4$. But we know $\overline{D} \in \Syl_2(E(\overline{G}))$ has rank $2$. This forces $\overline{G}$ to have a unique component $\overline{L}$. Hence, $\overline{D}$ is a Sylow $2$-subgroup of a simple group $E(\overline{G}) = \overline{L}$. Then Lemma \ref{lm:sylow-choices} gives us the result.
\end{proof}

\begin{lemma} \label{lm:unique-block}
    We have $O_{2'}(G) = 1$ and $C_G(O_2(G)) \leq O_2(G)$. In particular, $F^*(G) = O_2(G)$, and $B$ is the unique block of $G$.
\end{lemma}
\begin{proof}
    Since $F^*(\overline{G}) = E(\overline{G}) \cong \PSL_2(q)$ and $|\overline{G} : F^*(\overline{G})|$ has odd order, we deduce that $\overline{G} \cong \PSL_2(q)$. Then $\overline{G}$ has Schur Multiplier $C_2$. Since $B$ is reduced, it follows that $O_{2'}(G) \leq Z(G) \cap G'$. We deduce that $O_{2'}(G) = 1$. Thus, $O_2(G) = F^*(G)$, and we conclude that $C_G(O_2(G)) \leq O_2(G)$. This means that $B$ is the unique block of $G$ (see, for example, \cite[Theorem 4.14]{book:navarro}).
\end{proof}

\begin{lemma} \label{lm:resistant-no-comp}
    Let $\calF := \calF_D(G)$. Then $O_2(\calF) = O_2(G) = \Phi(D)$. In particular, $D \not\normalIn \calF$.
\end{lemma}
\begin{proof}
    We have found that $E(\overline{G}) \cong \PSL_2(q)$ for some $q \equiv \pm 3 \bmod{8}$ and that $F^*(G) = O_2(G)$. Let $\calF := \calF_D(G)$, which satisfies $O_2(G) \normalIn \calF$, and set $A := \Phi(D)$. We note that a Sylow $2$-subgroup of $\overline{G}$ is elementary abelian of order $2^2$, meaning that $A = O_2(G)$. In particular, $C_G(A) \leq A$. We deduce that $\Out_\calF(A) = \Out_G(A) \cong \PSL_2(q)$. Thus, $D \not\normalIn \calF$.
\end{proof}

\begin{proof}[Proof of Proposition \ref{prp:no-comp}]
    We have found in Lemma \ref{lm:unique-block} that $B$ is the principal block of $G$. As such, Lemma \ref{lm:resistant-no-comp} implies that $B$ cannot be controlled, which is a contradiction.
\end{proof}

\section{Components} \label{sec:comp}
In this section, we will prove the following result.

\begin{proposition} \label{prp:min-counter}
    Let $D$ be a maxperm $2$-group of rank $2$ such that $|D| \geq 2^4$. If $G$ is a quasisimple group and $B$ is a non-nilpotent, controlled block of $G$, then $D$ cannot be isomorphic to a defect group of $G$.
\end{proposition}
\noindent We make use of the work in Section \ref{sec:quasi-ab} to conclude this result.

We remark that if $|D| \leq 2^3$, then we have possible maxperm groups that are defect groups of quasisimple groups. In particular, a defect group of the principal block of $\PSL_2(q)$ (for $q \equiv \pm 3 \bmod{8}$) is elementary abelian of order $4$. Moreover, the principal block of the central extension $\SL_2(q)$ (for $q \equiv \pm 3 \bmod{8}$) has defect group $Q_8$. 

We argue by contradiction. Suppose that a quasisimple group $G$ has a non-nilpotent, controlled block $B$ with a defect group $D$ that is a maxperm $2$-group. Like in the previous section, we prove this in multiple steps. We write $\overline{G} := G/Z(G)$. Our proof loosely follows \cite[Proposition 8]{paper:ekb-block-mna}.
\begin{lemma}
    $D$ is not a Sylow $2$-subgroup of $G$. In particular, $|G|_2 \geq 2^5$.
\end{lemma}
\begin{proof}
    Since $D$ is maxperm, $N_G(D)/C_G(D)$ is not a $2$-group. In particular, there exists a $3$-element $\theta \in \Aut_G(D)$. Since $D$ has rank $2$, it follows that $\overline{D}$ is either cyclic or has rank $2$. We know that a Sylow $2$-subgroup of a simple group cannot be cyclic. Thus, $\overline{D}$ must have rank $2$. As $\overline{G}$ is simple, $O_2(Z(G))$ cannot be a Sylow $2$-subgroup of $G$. Since $O_2(Z(G))$ is invariant under $N_G(D)$, we deduce by Lemma \ref{lm:maxperm-char-in-frat} that $O_2(Z(G)) \leq \Phi(G)$. Then Lemma \ref{lm:burnside} implies that $\Aut(\overline{D})$ is not a $2$-group, and so $\overline{D}$ is a maxperm $2$-group by Lemma \ref{lm:maxperm-rank2}. We now appeal to Lemma \ref{lm:sylow-choices} to deduce that $\overline{G} \cong \PSL_2(q)$ for $q \equiv \pm 3 \bmod{8}$. Then $|Z(G)| \leq 2$, which implies that $|D| \leq 2^3$. This contradicts our assumption that $|D| \geq 2^4$.
\end{proof}
    
\begin{lemma}
    $\overline{G}$ is not an alternating group.
\end{lemma}
\begin{proof}
    Write $\overline{G} \cong A_n$. Since $|G|_2 \geq 2^5$, it follows that $n \geq 8$. As such, $|Z(G)| \leq 2$. 
    
    By Proposition \ref{prp:block-domination}, there exists a block $\overline{B}$ of $\overline{G}$ with defect group $\overline{D}$. We write $\overline{G_0} \cong S_n$ so that $O^2(\overline{G_0}) = \overline{G}$. Then there exists a block $\overline{B_0}$ of $\overline{G_0}$ that covers $\overline{B}$. Let $\overline{D_0}$ be a defect group of $\overline{B_0}$ satisfying $\overline{D} = \overline{D_0} \cap \overline{G}$. If $\overline{B_0}$ has weight $w$, then \cite{book:symmetric} reveals that $\overline{D_0}$ is $\overline{G_0}$-conjugate to a Sylow $2$-subgroup of $\overline{A_0} \cong S_{2w}$. If $\overline{A} := O^2(\overline{A_0}) \cong A_{2w}$, then it follows that
    \[\overline{D} = \overline{D_0} \cap \overline{G} = \overline{D_0} \cap \overline{A_0} \cap \overline{G} = \overline{D_0} \cap \overline{G_0}.\]
    This means that $\overline{D}$ is a Sylow $2$-subgroup of $\overline{A_0}$. Since $D$ is maxperm, we deduce that $\overline{D}$ must also be a maxperm $2$-group of rank $2$ by Lemma \ref{lm:maxperm-quo}. We infer that $4 \leq w \leq 5$. But then $|D| \leq 2|\overline{D}| \leq 2^3$, which is a contradiction.
\end{proof}

\begin{lemma}
    $\overline{G}$ is not a sporadic group.
\end{lemma}
\begin{proof}
    Using \cite[Table 2]{paper:nilpotent-abelian-2} and the proof of \cite[Proposition 4.1]{paper:nilpotent-abelian-2}, we see that $\overline{G}$ must be isomorphic to $\Ly$ or $\HN$. In particular, we have that $|Z(G)| = 1$. If $G \cong \HN$, then we find that $D \cong \SD_{16}$ by \cite[Table 10]{paper:awc}. Hence, $D$ is not a maxperm $2$-group, a contradiction. On the other hand, if $G \cong \Ly$, then $|D| = 2^6$. However, we compute in GAP that $D$ has rank $3$ in this case, contradicting our hypothesis again.
\end{proof}

\begin{lemma}
    $\overline{G}$ is not a group of Lie type in characteristic $2$.
\end{lemma}
\begin{proof}
    By \cite{paper:def-char}, we know that a $2$-block $B$ of $G$ has either full or zero defect. If $B$ has zero defect, it is nilpotent. On the other hand, if $B$ has full defect, then a defect group is a Sylow $2$-subgroup of $G$, a contradiction.
\end{proof}

We have found that if a maxperm $2$-group of rank $2$ $D$ is a defect group of a non-nilpotent block $B$ quasisimple group $G$, then $G$ must be a group of Lie type in odd characteristic. We will rule this particular case out in multiple steps.

\begin{lemma} \label{lm:smaller-sl23}
    Let $H$ be a finite group with $Z \leq Z(H)$ a $2$-group, and write $\overline{H} := H/Z$. Let $M$ and $T$ be subgroups of $H$ containing $Z$ such that:
    \begin{enumerate}
        \item $\overline{H} = \overline{M} \overline{T}$;
        \item $\overline{M} \normalIn \overline{H}$;
        \item $\overline{M}$ is a central product of groups of Lie type $\overline{M}_i$, for $1 \leq i \leq n$, in odd characteristic; and
        \item $\overline{T}$ is an abelian $2$-group.
    \end{enumerate}
    If there exists an $1 \leq i \leq n$ such that $\overline{M_i}/Z(\overline{M_i})$ is isomorphic to $\PSL_2(3)$, then any block $B$ of $H$ cannot be nilpotent.
\end{lemma}
\begin{proof}
    Assume that $B$ is nilpotent. By Proposition \ref{prp:block-domination}, there exists a non-nilpotent block $\overline{B}$ of $\overline{H}$ with defect group $\overline{D}$. Using (4), we infer that $\overline{G}/\overline{M}$ is a $2$-group. Applying Lemma \ref{lm:idxp-nilpotent}, we deduce that the block $\overline{b}$ of $\overline{M}$ covered by $\overline{B}$ is not nilpotent. Moreover, by Lemma \ref{lm:cen-prod-nilpotent}, it follows that the block $\overline{b}_i$ of $\overline{M_i}$ covered by $\overline{b}$ is nilpotent. But $\overline{M_i}/Z(\overline{M_i})$ is isomorphic to $\PSL_2(3)$. This implies that the block $\overline{b_i}$ is the unique block of $\overline{M_i}$. By Brauer's Third Main Theorem, it follows that $b_i$ is not nilpotent, which is a contradiction.
\end{proof}

\begin{lemma} \label{lm:smaller-non-nilpotent}
    Let $H$ be a finite group with $Z \leq Z(H)$ a $2$-group, and write $\overline{H} := H/Z$. Let $M$ and $T$ be subgroups of $H$ containing $Z$ such that:
    \begin{enumerate}
        \item $\overline{H} = \overline{M} \overline{T}$;
        \item $\overline{M} \normalIn \overline{H}$;
        \item $\overline{M}$ is a central product of groups of Lie type $\overline{M}_i$, for $1 \leq i \leq n$, in odd characteristic; and
        \item $\overline{T}$ is an abelian $2$-group.
    \end{enumerate}
    Assume that $D$ is a defect group of a non-nilpotent block $B$ of $H$ and is not quasi-abelian. Then there exists a quasisimple group $M_0 \leq H$ with a non-nilpotent block $b_0$ that has $D$ as a defect group.
\end{lemma}
\begin{proof}
    Since $Z$ is a $2$-group and (4) holds, it follows that $H/M$ is a $2$-group. By Lemma \ref{lm:idxp-nilpotent}, we deduce that the block $b$ of $M$ covered by $B$ is not nilpotent. Applying Lemma \ref{lm:smaller-not-nilpotent}, we may assume that $D$ is a defect group of $b$.

    Using Proposition \ref{prp:block-domination}, there exists a non-nilpotent block $\overline{b}$ of $\overline{M}$ with defect group $\overline{D}$. Then Lemma \ref{lm:cen-prod-nilpotent} reveals that there is a non-nilpotent block $\overline{b}_i$ of some $\overline{M_i}$ covered by $\overline{b}$. This corresponds to a non-nilpotent block $b_i$ of $M_i$, where $M_i$ is the full preimage of $\overline{M_i}$ in $H$. By construction, $b_i$ covers $M$. Assume further that $D \cap M_i$ is a defect group of $b_i$. Applying Lemma \ref{lm:smaller-not-nilpotent}, we deduce that $D \leq M_i$.

    As $Z \leq O_2(Z(H))$, we infer that $Z \leq Z(D)$. If $\overline{M_i}/Z(\overline{M_i}) \cong \PSL_2(3)$, then it follows that $|D : Z| \leq 2^3$. By Lemma \ref{lm:maxperm-small}, we infer that $D$ is minimal non-abelian. But then $D$ is quasi-abelian by Lemma \ref{lm:minab-qab}, which is a contradiction. Thus, $\overline{M_i}$ must be a quasisimple group.

    Set $M_0 := M_i'$. Let $b_0$ be a block of $M_0$ covered by $b_i$. Since $Z$ is a $2$-group, Lemma \ref{lm:perfect-cen} reveals that $M_i/M_0$ is a $2$-group. This implies that $b_0$ is not nilpotent by Lemma \ref{lm:idxp-nilpotent}. Applying Lemma \ref{lm:smaller-not-nilpotent}, we deduce that $D$ is contained inside $M_0$. We further know that $M_0$ is a perfect group by Lemma \ref{lm:perfect-cen}, making $M_0$ quasisimple.
\end{proof}

\begin{lemma} \label{lm:smaller-nilpotent}
    Let $H$ be a finite group with subgroups $M$ and $T_0$ such that:
    \begin{enumerate}
        \item $H = MT_0$;
        \item $M \normalIn H$;
        \item $M$ is a central product of groups $M_i$ of quasisimple groups; and
        \item $T_0$ a $2$-group such that $T_0 Z(H)/Z(H)$ is abelian.
    \end{enumerate}
    If $D$ is a defect group of a nilpotent block $B$ of $H$, then $D$ is quasi-abelian.
\end{lemma}
\begin{proof}
    Let $b$ be a block of $M$ covered by $B$. We may assume that a defect group of $b$ is $D \cap M$. We set $b_i$ to be a block of $M_i$, and assume that a defect group of $b_i$ is $D \cap M_i$.
    
    Applying Lemma \ref{lm:idxp-nilpotent}, we deduce that the block $b$ of $M$ is also nilpotent. By Lemma \ref{lm:cen-prod-nilpotent}, each $b_i$ is also nilpotent. But each $M_i$ is a quasisimple group, so \cite{paper:nilpotent-abelian-2} tells us that $D \cap M_i$ is abelian for $1 \leq i \leq n$. We recall that $M$ is a central product of $M_i$. This implies that $A := D \cap M$ is abelian as well. As such, $D = AT$, where $T := T_0 \cap D$. We see that 
    \[\frac{T}{T \cap Z(D)} \cong \frac{T/(T \cap Z(H)))}{(T \cap Z(D))/(T \cap Z(H))}\] 
    is abelian, meaning that $D$ is quasi-abelian by Lemma \ref{lm:quasiab-equiv}.
\end{proof}

We recall that $G$ is a quasisimple group with a non-nilpotent block $B$ such that a non-abelian maxperm $2$-group $D$ of rank $2$ is a defect group of $B$. Using our analysis in this section, $\overline{G}$ must be a group of Lie type in odd characteristic. Moreover, $\overline{G}$ is not isomorphic to an alternating group or a group of Lie type in even characteristic. 

\begin{lemma}
    $\overline{G}$ is not a group of Lie type in odd characteristic.
\end{lemma}
\begin{proof}
    We may assume that $|G|$ is minimal amongst all quasisimple groups with a non-nilpotent block $B$ that has $D$ as a defect group. We first show that $D$ is quasi-abelian. 
    
    Assume that $G/Z(G) \cong \PSU_4(3)$. We check using GAP \cite{GAP4} that any $2$-subgroup $P$ of $G$ that is maxperm of rank $2$ satisfies $|P| \leq 2^6$. As such, we have $|D| \leq 2^6$, meaning that $D$ is quasi-abelian by Lemma \ref{lm:qab-small}. For the remaining groups, $G/Z(G)$ does not have exceptional cover. Hence, $G$ is a group of Lie type in odd characteristic.
    
    Let $Z := O_2(Z(G))$ and write $\overline{G} := G/Z$. Let $\overline{P} \in \Syl_2(\overline{G})$ containing $\overline{D}$, and fix $t \in Z(\overline{P})$. Let $Z \leq H \leq G$ be such that $\overline{H} := C_{\overline{G}}(t)$. Since $t$ is not central in $\overline{G}$, $\overline{H}$ is a proper subgroup of $\overline{G}$ containing $\overline{D}$. By \cite[Theorem 4.2.2]{book:gls3}, there exist subgroups $M, T_0$ of $G$ containing $Z$ such that $\overline{H} = \overline{M} \overline{T_0}$, where $\overline{M}$ is a central product of groups of Lie type in odd characteristic, and $\overline{T_0}$ is an abelian group. By replacing $\overline{T_0}$ with $O_2(\overline{T_0})$, we may further assume that $\overline{T_0}$ (and $T_0$) is a $2$-group. 

    Let $b$ be a block of $H = MT_0$ with defect group $D$. We note that $|G|$ was chosen to be minimal amongst all quasisimple groups $G_0$ with a non-nilpotent that has $D$ as a defect group. Thus, we may apply Lemma \ref{lm:smaller-non-nilpotent} to deduce that either $D$ is quasi-abelian or that the block $b$ is nilpotent. Assume that the block $b$ is nilpotent. Then Lemma \ref{lm:smaller-sl23} implies that $\overline{M_i}/Z(\overline{M_i})$ is simple for each $1 \leq i \leq n$. In particular, $\overline{M}$ is a perfect group. By Lemma \ref{lm:perfect-cen}, we deduce that $M_0 := M'$ is a perfect group and that $M = M_0 Z$. Furthermore, Lemma \ref{lm:grun} implies that $M_0/Z(M_0)$ is a direct product of simple groups. As such, $M_0$ is a central product of quasisimple groups. Moreover, since $Z \leq T_0$, we get that 
    \[M_0 T_0 = M_0 ZT_0 = MT_0 = H.\]
    We now assume that $M = M_0$. Since $b$ is nilpotent, we deduce by Lemma \ref{lm:smaller-nilpotent} that $D$ is quasi-abelian. Thus, $D$ is quasi-abelian in every case.

    Since $D$ is also a maxperm $2$-group of rank $2$, it follows that $D$ is resistant by Proposition \ref{prp:metabelian-res}. This implies that the block $B$ with defect group $D$ is controlled. Since $|D| \geq 2^4$, we get a contradiction by \cite[Theorem 1.1]{paper:j-an-controlled-2block}.
\end{proof}

\section{Main Proofs} \label{sec:final}
We are now in a position to prove the main results.

\begin{theorem} \label{thm:main}
    Let $G$ be a finite group, and let $B$ be a non-nilpotent, controlled block of $G$ with defect group $D$, where $|D| \geq 2^4$ is a maxperm $2$-group of rank $2$. Then $B$ is Morita equivalent to $\calO [D \rtimes C_3]$.
\end{theorem}
\begin{proof}
    If $D$ is abelian, then we appeal to Lemma \ref{lm:maxperm-quo} to see that $D$ is homocyclic group of rank $2$. In that case, \cite{paper:hom-rank-2} reveals that $B$ is Morita equivalent to $\calO[D \rtimes C_3]$. 
    
    We assume now that $D$ is non-abelian. By Proposition \ref{prp:reduced}, we may further assume that $(G, B)$ is reduced. If $O_2(G) = D$, then $N_G(D)/DC_G(D)$ is isomorphic to $C_3$ as $D$ has rank $2$. It follows that $B$ is Morita equivalent to $\calO[D \rtimes C_3]$ by \cite{paper:kul-normals}. 
    
    We shall now assume that $O_2(G) < D$. By Proposition \ref{prp:no-comp}, it follows that $E(G) \neq 1$. Let $B_E$ be a block of $E(G)$ covered by $B$. If $D \nleq E(G)$, then $B_E$ is nilpotent by Lemma \ref{lm:smaller-not-nilpotent}, which contradicts that $B$ is reduced. Hence, $B_E$ cannot be nilpotent. By Lemma \ref{lm:cen-prod-nilpotent}, we deduce that there exists a non-nilpotent block $B_L$ of a quasisimple group $L \normalIn E(G)$ covered by $B_E$. Applying Lemma \ref{lm:smaller-not-nilpotent} on $E(G)$, we infer that $D$ is a defect group of $L$. But then Proposition \ref{prp:min-counter} gives us a contradiction.
\end{proof}

\begin{corollary}
    Let $G$ be a finite group, and let $B$ be a block of $G$ with defect group $D$, a resistant maxperm $2$-group of rank $2$ with $|D| \geq 2^4$. Then $B$ is Morita equivalent to $\calO D$ or $\calO[D \rtimes C_3]$.
\end{corollary}
\begin{proof}
    Since $D$ is resistant, $B$ is a controlled block. If $B$ is nilpotent, then $B$ is Morita equivalent to $\calO D$ by \cite{paper:nilpotent}. Otherwise, Theorem \ref{thm:main} tells us that $B$ is Morita equivalent to $\calO[D \rtimes C_3]$.
\end{proof}
We recall that if $D$ is a $2$-generated maxperm $2$-group satisfying $|D| \leq 2^3$, then $D$ is isomorphic to $C_2 \times C_2$ or $Q_8$. If $B$ is a block of some finite group with defect group $D \cong C_2 \times C_2$, then \cite{paper:li-c2c2} reveals that $B$ is Morita equivalent to $\calO D$, $\calO A_4$ or the principal block of $\calO A_5$. If $B$ is a block with defect group $D \cong Q_8$, then \cite{paper:ei-q8} implies that $B$ is Morita equivalent to $\calO D$, $\calO \SL_2(3)$ or the principal block of $\calO \SL_2(5)$.

\section*{Acknowledgments}
This work is part of the author's PhD at the University of Manchester under the supervision of Prof. Charles Eaton and Dr. Martin van Beek. The author would also like to thank Prof. Radha Kessar for helpful discussions. The author acknowledges the support received from EPSRC (EP/W524347/1) during this period. 

\bibliographystyle{plain}
\bibliography{nat}
\end{document}